\documentclass[12pt,twoside]{amsart}
\usepackage{amssymb,enumerate}

\numberwithin{equation}{section}
\newtheorem{theorem}{Theorem}[section]
\newtheorem{lemma}[theorem]{Lemma}
\newtheorem{proposition}[theorem]{Proposition}
\newtheorem{corollary}[theorem]{Corollary}

\newenvironment{theoremnonumber}{\begin{trivlist}\item[]{\bf Theorem.}\ \em}
{\end{trivlist}}

\theoremstyle{definition}
\newtheorem{definition}[theorem]{Definition} 
\newtheorem{remark}[theorem]{Remark}
\newtheorem{example}[theorem]{Example}

\newtheorem{question}[theorem]{Question}

\begin{document}


\newcommand{\m}[1]{\marginpar{\addtolength{\baselineskip}{-3pt}{\footnotesize \it #1}}}
\newcommand{\A}{\mathcal{A}}
\newcommand{\K}{\mathcal{K}} 
\newcommand{\knd}{\mathcal{K}^{[d]}_n}
\newcommand{\F}{\mathcal{F}}
\newcommand{\N}{\mathbb{N}}
\newcommand{\pr}{\mathbb{P}}
\newcommand{\I}{\mathcal{I}}
\newcommand{\G}{\mathcal{G}}
\newcommand{\lcm}{\operatorname{lcm}}
\newcommand{\ndp}{N_{d,p}}
\newcommand{\tor}{\operatorname{Tor}}
\newcommand{\reg}{\operatorname{reg}} 
\newcommand{\mf}{\mathfrak{m}}
\def\bb{{{\rm \bf b}}}
\def\cc{{{\rm \bf c}}}

 
\title{Tetrahedral curves via graphs and Alexander duality}
\author{Christopher A. Francisco}
\address{Department of Mathematics, University of Missouri, Mathematical Sciences Building, Columbia, MO 65203}
\email{chrisf@math.missouri.edu}
\urladdr{http://www.math.missouri.edu/$\sim$chrisf}


\begin{abstract}
A tetrahedral curve is a (usually nonreduced) curve in $\pr^3$ defined by an unmixed, height two ideal generated by monomials. We characterize when these curves are arithmetically Cohen-Macaulay by associating a graph to each curve and, using results from combinatorial commutative algebra and Alexander duality, relating the structure of the complementary graph to the Cohen-Macaulay property.
\end{abstract}
 
\maketitle

\section{Introduction} \label{s.intro}

Our purpose in this paper is to demonstrate the use of several techniques from combinatorial commutative algebra in order to answer a question from \cite{MN} about algebraic properties of certain curves. The curves we study lie in $\pr^3$ and are defined by ideals of the form \[ I = (a,b)^{p_1} \cap (a,c)^{p_2} \cap (a,d)^{p_3} \cap (b,c)^{p_4} \cap (b,d)^{p_5} \cap (c,d)^{p_6} \subset k[a,b,c,d], \] where $k$ is a field, and the $p_i$ are nonnegative integers. These curves are called {\it tetrahedral curves} because one can view the six lines defined by the height two ideals as the edges of a tetrahedron. In this work, we characterize explicitly when these curves are arithmetically Cohen-Macaulay (meaning that their coordinate ring is Cohen-Macaulay, and abbreviated throughout as ACM) by repeatedly reformulating this property for tetrahedral curves in different combinatorial and numerical terms.

The study of tetrahedral curves dates back at least to Schwartau's 1982 Ph.D. thesis \cite{S} in which he investigated the case $p_2=p_5=0$. Schwartau's main interest was to determine when such tetrahedral curves are ACM with the technique of liaison addition. Recently, Migliore and Nagel generalized this work substantially. Using basic double linkage, a special case of liaison addition, Migliore and Nagel carried out a comprehensive study of the properties of tetrahedral curves in \cite{MN}. They developed a reduction procedure for tetrahedral curves to investigate a number of questions, including when such curves are ACM, the structure of their minimal free resolutions, and some properties of the Hilbert scheme. The author, Migliore, and Nagel, extended this work in \cite{FMN} to characterize when the ideal of a tetrahedral curve is componentwise linear, give a numerical algorithm to determine its graded Betti numbers, and partially describe its generic initial ideal. 

Migliore and Nagel's reduction algorithm in \cite{MN} is an efficient mechanism for determining when a particular tetrahedral curve is arithmetically Cohen-Macaulay. Given a tetrahedral curve with exponent vector $(p_1,\dots,p_6)$, one carries out a series of reductions of the $p_i$ according to rules governed by taking basic double links. Eventually, one either reaches the {\it trivial curve} defined by the exponent vector $(0,\dots,0)$, or one reaches a {\it minimal curve} that is not the trivial curve, and one cannot do any further reductions. In the first case, the original curve is ACM; in the second, it is not. The variety of applications of Migliore and Nagel's algorithm in \cite{MN} is an excellent illustration of the effectiveness of the basic double linkage technique.

The tetrahedral curves with which Schwartau worked, which have $p_2=p_5=0$ and defining ideals of the form \[ I=(a,b)^{p_1} \cap (a,d)^{p_3} \cap (b,c)^{p_4} \cap (c,d)^{p_6},\] are called {\it Schwartau curves}. Both Schwartau \cite{S} (quoted in \cite[Theorem 2.4]{MN}) and Migliore and Nagel \cite[Theorem 5.3]{MN} gave explicit necessary and sufficient conditions on $p_1,p_3,p_4$, and $p_6$ for the Schwartau curve to be ACM. In \cite{MN}, this result is a consequence of a careful analysis of Migliore and Nagel's algorithm, computing which exponent vectors reduce to $(0,\dots,0)$.

In \cite[Question 7.4(5)]{MN}, Migliore and Nagel ask about the natural generalization of these results:

\begin{question} \label{q.mn}
Can the tetrahedral curves in $\pr^3$ that are arithmetically Cohen-Macaulay be identified by explicitly giving the 6-tuples (as Schwartau does for 4-tuples [in the case of Schwartau curves])?
\end{question}

That is, can we find explicit conditions on the $p_i$ (other than saying that the vector $(p_1,\dots,p_6)$ reduces to all zeros in the Migliore-Nagel algorithm) that will tell us precisely when a tetrahedral curve is ACM? We would like simply to have some inequalities on the $p_i$ to test that would tell us immediately whether or not a curve is ACM, avoiding the need to use any numerical algorithm. This is a substantially more difficult question than in the Schwartau curve case even with Migliore and Nagel's algorithm available.

In this paper, we use a very different approach and give a positive answer to Question~\ref{q.mn}. Our main theorem is Theorem~\ref{t.main}, which generalizes Schwartau's work:

\begin{theoremnonumber}
Let $C$ be a tetrahedral curve with exponent vector $(p_1,\dots,p_6)$. Suppose without loss of generality that $p_1+p_6=\max(p_1+p_6,p_2+p_5,p_3+p_4)$. Then $C$ is ACM if and only if at least one of the following conditions holds:
\begin{enumerate}[(i)]
\item $p_1=0$ or $p_6=0$ 
\item $p_1+p_6=\epsilon+\max(p_2+p_5,p_3+p_4)$, where $\epsilon \in \{0,1\}$
\item $2p_1 < p_2+p_3+3-p_6$ or $2p_1 < p_4+p_5+3-p_6$ or $2p_6 < p_2+p_4+3-p_1$ or $2p_6 < p_3+p_5+3-p_1$
\item All inequalities of (iii) fail, $p_1+p_6=2+p_2+p_5=2+p_3+p_4$, and $p_1+p_3+p_5$ is even.
\end{enumerate}
\end{theoremnonumber}

While we think it is interesting by itself to have a complete answer to Migliore and Nagel's Question~\ref{q.mn}, even if two of the conditions are relatively unwieldy, we believe that the primary interest in this work is in the method we used to solve the problem, particularly because the necessary and sufficient conditions are not at all easy to isolate. Our approach to Question~\ref{q.mn} avoids liaison theory entirely; the techniques are combinatorial. We give a brief outline here of the methods used. For a more detailed discussion, see Section~\ref{s.equivalent}. 

Let $I$ be the ideal of a tetrahedral curve $C$ with exponent vector $(p_1,\dots,p_6)$. We proceed in the following manner:
\begin{itemize}
\item Polarize $I$ to get a squarefree monomial ideal $J$.
\item Take the Alexander dual of $J$, giving a new squarefree monomial ideal $J^{\vee}$.
\item Because $I$ and $J$ are unmixed of height two, $J^{\vee}$ is generated by degree two squarefree monomials. Hence it is the edge ideal of a graph $G$.
\item Given $G$, take the complementary graph $\bar{G}$, the graph on the same vertex set whose edges are precisely those not appearing in $G$.
\item Analyze what kind of induced $r$-cycles can appear in $\bar{G}$ for $r \ge 4$.
\item Determine numerical conditions in terms of the $p_i$ for $\bar{G}$ to be a chordal graph.
\item Combine these efforts with theorems of Eagon and Reiner from Alexander duality and Fr\"oberg on the resolutions of edge ideals to reformulate Question~\ref{q.mn} as a purely numerical question about the partitioning of two positive integers, subject to four inequalities.
\end{itemize}

Our method reframes Question~\ref{q.mn} as a problem about the structure of graphs associated to tetrahedral curves, which can be further interpreted as a set of numerical conditions on $p_1,\dots,p_6$. While the original problem of determining when a tetrahedral curve is ACM is difficult, the numerical translation obtained through combinatorial commutative algebra is completely elementary and requires only a detailed analysis of when two positive integers can be partitioned in a particular way. The reformulations also make it evident that for a tetrahedral curve, the ACM property is independent of the characteristic of the underlying field.

Our approach relies on polarization preserving Cohen-Macaulayness (or lack of it), Alexander duality's ability to translate between the Cohen-Macaulay property and the structure of the resolution of the Alexander dual, and the knowledge that the resolution of the edge ideal of a graph is greatly influenced by the structure of the complementary graph. Though these techniques have been used separately many times, we believe the combination of all of them in this way is new, and, in particular, gives a nice application of Fr\"oberg's work on graphs.

The method of investigation used here could be easily extended to the analogous unmixed, height two ideals in more variables. We note that the techniques in this paper depend on the fact that we know a reasonable amount about the resolutions of edge ideals, and thus these results present ideas for future applications of the work being done on resolutions of facet ideals of higher-dimensional simplicial complexes (see, e.g., \cite{Faridi-facet,HVT,HVTsurvey,Zheng}). A principal goal in this research program is to find generalizations of Fr\"oberg's work connecting graphs and free resolutions. Such advances would allow us to use the approach in this paper to investigate monomial ideals that are not unmixed or that have height greater than two; moreover, this gives further justification for studying properties of facet ideals of simplicial complexes in addition to their Stanley-Reisner ideals. We also refer the reader to recent papers such as \cite{HHZ,Katzman} that study graphs in the context of commutative algebra and resolutions and Villarreal's book \cite{V} for a good introduction to these ideas.

We organize the paper in the following manner. In Section~\ref{s.prelim}, we gather preliminary results about resolutions of edge ideals, polarization, and Alexander duality. We give some equivalent formulations of the statement that a tetrahedral curve is ACM in Section~\ref{s.equivalent}, reducing Question~\ref{q.mn} to a numerical question about partitioning integers. Using this reformulation, Section~\ref{s.acm} gives some sufficient conditions for a tetrahedral curve to be ACM, and we conclude in Section~\ref{s.final} by analyzing the remaining cases to complete the characterization of ACM tetrahedral curves.

We thank Juan Migliore and Uwe Nagel for introducing us to tetrahedral curves, many helpful conversations on the subject, and for posing the question that inspired this work. We gratefully acknowledge Grayson and Stillman's computer algebra system Macaulay 2 \cite{M2} that guided us in determining the conditions of Theorem~\ref{t.main} and suggested methods of proof. Finally, we thank the referee for his or her helpful comments and suggestions.

\section{Preliminaries on graphs, polarization, and Alexander duality} \label{s.prelim}

In this section, we state some definitions from graph theory and fundamental results on the resolutions of edge ideals and Alexander duality. These tools will allow us to make a sequence of reformulations of the question of when a tetrahedral curve is ACM.

\subsection{Some graph terminology} \label{sub.graph}

Let $G$ be a graph on vertex set $V=\{x_1,\dots,x_n\}$ with no loops or multiple edges, and let $E$ be the edge set of $G$, which is comprised of two-element subsets of $V$. We will be interested in particular types of subgraphs of $G$.

\begin{definition} \label{d.induced}
Let $G$ be a graph with vertex set $V$ and edge set $E$. A subgraph $G'$ of $G$ with vertex set $V' \subseteq V$ and edge set $E' \subseteq E$ is an \textbf{induced subgraph} of $G$ if for all $x_i$ and $x_j$ in $V'$, $\{x_i,x_j\} \in E'$ if and only if $\{x_i,x_j\} \in E$. The induced subgraph on $V' \subseteq V$ is the induced subgraph of $G$ with vertex set $V'$ and edge set $E' \subseteq E$, which consists of all edges of $G$ involving only vertices of $V'$.
\end{definition}

We use the notion of an induced subgraph to define a chordal graph. For notational ease, we often write an edge $\{x_i,x_j\}$ as $x_ix_j$. Recall that a cycle inside a graph $G$ is a sequence of distinct vertices $x_1, \dots, x_r$ in the vertex set of $G$ such that there exist edges $x_1x_2$, $x_2x_3, \dots$, $x_{r-1}x_r, x_1x_r$ in the edge set of $G$. 

\begin{definition} \label{d.chordal}
A graph $G$ is {\bf chordal} if the only induced cycles in $G$ are triangles. That is, for $r \ge 4$, any $r$-cycle inside $G$ has a chord, an edge connecting two nonconsecutive vertices in the cycle.
\end{definition}

Given a graph $G$, there exists a complementary graph whose structure is vital in studying the algebraic properties of the ideal associated to $G$.

\begin{definition} \label{d.complementary}
Let $G$ be a graph with vertex set $V$ and edge set $E$. The \textbf{complementary graph} of $G$ is a graph $\bar{G}$ with vertex set $V$. Its edges are precisely the complement of the edges of $G$; two vertices are connected by an edge in $\bar{G}$ if an only if there is no edge between them in $G$.
\end{definition}

In order to connect commutative algebra to the structure of a graph, we define the edge ideal of a graph.

\begin{definition} \label{d.edgeideal}
The edge ideal of a graph $G$ on vertices $x_1,\dots,x_n$ is the ideal \[ \I(G) = (x_ix_j : x_i \mbox{ and } x_j \, \mbox{ are connected by an edge in } G) \subset k[x_1,\dots,x_n].\] 
\end{definition}

Note that this is a facet ideal (omitting isolated vertices); it is not the Stanley-Reisner ideal of $G$ considered as a one-dimensional simplicial complex. There have been many papers relating algebraic properties of the edge ideal to concepts in graph theory. One particularly useful result for us is a theorem of Fr\"oberg from \cite{Froberg}.

\begin{theorem} \label{t.froberg}
Let $G$ be a graph with edge ideal $\I(G)$. Then $\I(G)$ has a linear resolution if and only if $\bar{G}$ is chordal.
\end{theorem}

Fr\"oberg's theorem is one of the two fundamental results that will allow us to use graphs to recharacterize Question~\ref{q.mn} as a purely numerical problem. Note that the characteristic of the underlying field $k$ is irrelevant in Theorem~\ref{t.froberg}. Finding an analogous result for higher-dimensional simplicial complexes will be harder because, for example, there are squarefree monomial ideals generated in degree three that have linear resolutions in characteristic not two but nonlinear minimal resolutions if char $k$ = 2 (see, e.g., the example of a triangulation of the real projective plane, due to Reisner, in \cite[page 228]{BH}). Generalizations of Theorem~\ref{t.froberg} will likely have to restrict to a fixed characteristic or to particular simplicial complexes whose homology is independent of the characteristic of the underlying field. 

\subsection{Polarization} \label{sub.polar}

We wish to use squarefree Alexander duality, a tool for understanding simplicial complexes (and, equivalently, squarefree monomial ideals), to investigate tetrahedral curves. The difficulty is that unless all $p_i$ are zero or one, the ideal defining a tetrahedral curve is not squarefree. In order to associate a squarefree monomial ideal to a tetrahedral curve, we polarize the defining ideal of the curve. For details on polarization, see Faridi's paper \cite{Faridi}; we sketch a few main points here.

Polarizing maps an arbitrary monomial ideal to a squarefree monomial ideal in a polynomial ring with more variables. The operation changes a power of a variable into a product of many variables.

\begin{example} \label{e.polar}
Let $I=(x^3,x^2y^3,xy^4,y^5) \subset k[x,y]$. The polarization $J$ of $I$ is an ideal of the polynomial ring $k[x_1,x_2,x_3,y_1,\dots,y_5]$, and \[ J = (x_1x_2x_3, x_1x_2y_1y_2y_3,x_1y_1y_2y_3y_4, y_1y_2y_3y_4y_5).\]
\end{example}

Polarization has a number of helpful properties, including the following two results that we will use. See, for example, \cite[Propositions 2.8, 2.3(3),2.5(2)]{Faridi}.

\begin{lemma} \label{l.polar}
Let $R=k[x_1,\dots,x_n]$.
\begin{enumerate}[(i)]
\item Let $I$ be a monomial ideal in $R$, and let $J \subset R'$ be its polarization. Then $R/I$ is Cohen-Macaulay if and only if $R'/J$ is Cohen-Macaulay.
\item Let $M$ and $N$ be monomial ideals of $R$ with polarizations $M'$ and $N'$. Then the polarization of $M \cap N$ is equal to $M' \cap N'$; that is, polarization commutes with taking intersections.
\item If $I=(x_{i_1},\dots,x_{i_r})^m$ has polarization $J \subset R'$, then the associated primes of $R'/J$ are the ideals $(x_{i_1,c_1}, \dots, x_{i_r,c_r})$, where $1 \le c_j \le m$ for all $c_j$, and $c_1 + \cdots + c_r \le m+r-1$.
\end{enumerate}
\end{lemma}

\begin{example} \label{e.polar211112}
Let $C$ be the tetrahedral curve with exponent vector $(2,1,1,1,1,2)$. The ideal $I$ of $C$ is \[I=(a,b)^2 \cap (a,c) \cap (a,d) \cap (b,c) \cap (b,d) \cap (c,d)^2 = (abd^2, b^2cd, abcd, a^2cd, abc^2).\] Polarizing, we get the ideal \[ J = (a_1b_1d_1d_2, b_1b_2c_1d_1, a_1b_1c_1d_1, a_1a_2c_1d_1, a_1b_1c_1c_2).\] 
\end{example}

\subsection{Alexander duality} \label{sub.alex}

We conclude our section of preliminaries with a brief discussion of Alexander duality for squarefree monomial ideals. Given a squarefree monomial ideal $J$, we can compute its Alexander dual $J^{\vee}$ by mapping minimal generators of $J$ to components of $J^{\vee}$. 

\begin{definition} \label{d.alexdual}
If $J = (x_{1,1}x_{1,2}\cdots x_{1,t_1}, \ldots, x_{s,1}x_{s,2}\cdots x_{s,t_s})$ is a squarefree monomial ideal, then the {\bf Alexander dual} of $J$, denoted $J^{\vee}$, is the monomial ideal \[J^{\vee} = (x_{1,1},\ldots,x_{1,t_1}) \cap \cdots \cap (x_{s,1},\ldots,x_{s,t_s}).\]
\end{definition}

This definition is derived from Alexander duality on simplicial complexes: Given a simplicial complex $\Delta$, the Alexander dual of $\Delta$ is the simplicial complex \[\Delta^*=\{\{1,\dots,n\} \backslash F: F \not \in \Delta\}.\] If $\Delta^{*}$ is the Alexander dual of $\Delta$, then the Alexander dual of the Stanley-Reisner ideal $I_{\Delta}$ is the ideal $I_{\Delta}^{\vee} = I_{\Delta^*}$. 

Because the ideal of a tetrahedral curve is an unmixed, height two ideal, the polarization ideal $J$ is also unmixed of height two. Therefore, since Alexander duality maps generators to components, the Alexander dual of $J$ is generated by squarefree monomials of degree two.

\begin{example} \label{e.211112second}
Let $C$ be the tetrahedral curve $(2,1,1,1,1,2)$. In Example~\ref{e.polar211112}, we found that the ideal $I$ of $C$ has polarization \[ J = (a_1b_1d_1d_2, b_1b_2c_1d_1, a_1b_1c_1d_1, a_1a_2c_1d_1, a_1b_1c_1c_2) = \] \[(a_1,b_1) \cap (a_2,b_1) \cap (a_1,b_2) \cap (a_1,c_1) \cap (a_1,d_1) \cap (b_1,c_1) \cap (b_1,d_1) \cap (c_1,d_1) \cap (c_2,d_1) \cap (c_1,d_2).\] Mapping components to minimal generators, the Alexander dual of $J$ is \[ J^{\vee} = (a_1b_1, a_2b_1, a_1b_2, a_1c_1, a_1d_1, b_1c_1, b_1d_1, c_1d_1, c_2d_1, c_1d_2).\]
\end{example}

Alexander duality is an important tool in determining when quotients by squarefree monomial ideals (equivalently, simplicial complexes) are Cohen-Macaulay. The following theorem of Eagon and Reiner \cite{ER} allows us to translate between the Cohen-Macaulay property and the resolution of the Alexander dual.

\begin{theorem} \label{t.er}
Let $I$ be a squarefree monomial ideal in $R=k[x_1,\dots,x_n]$. Then $R/I$ is Cohen-Macaulay over $k$ if and only if $I^{\vee}$ has a linear resolution over $R$.
\end{theorem}

\begin{corollary} \label{c.linres}
Let $C$ be a tetrahedral curve with defining ideal $I$. Suppose $J$ is the polarization of $I$ and $J^{\vee}$ is the Alexander dual of $J$. Then $C$ is ACM if and only if $J^{\vee}$ has a linear resolution.
\end{corollary}

By Corollary~\ref{c.linres}, the curve with exponent vector $(2,1,1,1,1,2)$ of Example~\ref{e.211112second} is arithmetically Cohen-Macaulay because the Alexander dual $J^{\vee} \subset R'$ has minimal graded free resolution \[ 0 \longrightarrow R'(-5)^4 \longrightarrow R'(-4)^{15} \longrightarrow R'(-3)^{20} \longrightarrow R'(-2)^{10} \longrightarrow J^{\vee} \longrightarrow 0, \] which is linear.

\section{Equivalent formulations} \label{s.equivalent}

We reformulate the question of when a tetrahedral curve is ACM in several different ways in this section. Throughout, $C$ is a tetrahedral curve in $\pr^3$, defined by the ideal \[I=(a,b)^{p_1} \cap (a,c)^{p_2} \cap (a,d)^{p_3} \cap (b,c)^{p_4} \cap (b,d)^{p_5} \cap (c,d)^{p_6}.\] We will abuse notation, frequently writing $C=(p_1,\dots,p_6)$. In addition, our \textbf{standing assumption} throughout is that $p_1+p_6 \ge \max(p_2+p_5,p_3+p_4)$; if not, we can simply permute the variables to make the inequality true. The results of the last section show that we can start with $I$, polarize to get a squarefree ideal $J$, and take the Alexander dual to obtain an ideal $J^{\vee}$, which is the edge ideal of a graph $G$. We will refer to $G$ as the graph of $C$, and we will call the complementary graph $\bar{G}$ of $G$ the complementary graph of $C$.

Corollary~\ref{c.linres} translates the problem of determining when a tetrahedral curve is ACM to a question about when the resolution of a particular type of squarefree monomial ideal is linear. Combining Theorem~\ref{t.froberg} with Corollary~\ref{c.linres}, we immediately obtain:

\begin{proposition} \label{p.acmgraph}
A tetrahedral curve $C$ is ACM if and only if its complementary graph $\bar{G}$ is chordal.
\end{proposition}

\begin{remark} \label{r.morevars}
We remark that in the translations of the ACM property to the complemenatary graph characterization of Proposition~\ref{p.acmgraph}, there is nothing special about the fact that the ideal of a tetrahedral curve lives in $k[a,b,c,d]$. Identical arguments work for any unmixed, height two monomial ideal \[ I = \bigcap_{1 \le i < j \le n} (x_i,x_j)^{p_{i,j}} \subset k[x_1,\dots,x_n],\] where the $p_{i,j} \ge 0$. Consequently, Cohen-Macaulayness of these ideals is independent of the field $k$ because chordality of the associated complementary graph does not depend on $k$. 
\end{remark}

By Proposition~\ref{p.acmgraph}, to determine which tetrahedral curves are ACM, we need to characterize which complementary graphs $\bar{G}$ are chordal in terms of the exponent vector $(p_1,\dots,p_6)$. We spend the remainder of the section translating this into a purely numerical condition. (One could do a similar translation in more variables. We do not do this, however, because the numerical result would be much more complicated, and, we believe, would not shed further light on the techniques we wish to demonstrate.) The next result describes the minimal generating set of the Alexander dual of the polarization of the ideal of a tetrahedral curve, allowing us to get information about the graph associated to the curve and its complementary graph. 

\begin{lemma} \label{l.dualgens}
Let $I$ be the ideal of the tetrahedral curve $C=(p_1,\dots,p_6)$, and let $J^{\vee}$ be the Alexander dual of the polarization of $I$ in the variables $a_i$, $b_i$, $c_i$, and $d_i$. Then $a_ib_j$ is a minimal generator of $J^{\vee}$ if and only if $i$ and $j$ are positive integers with $i+j \le p_1+1$. The analogous characterization holds for the other types of generators.
\end{lemma}

\begin{proof}
To find the minimal generating set of $J^{\vee}$, it suffices to exhibit all associated primes of the polarization of $I$ since, for example, $a_ib_j$ is a minimal generator of $J^{\vee}$ if and only if $(a_i,b_j)$ is an associated prime of the polarization of $I$. The result thus follows from Lemma~\ref{l.polar}(ii), allowing us to compute the polarization of $I$ one component at a time, and Lemma~\ref{l.polar}(iii), which describes the primary decomposition of the polarization of an ideal of the form $(x_{i_1},\dots,x_{i_r})^m$; in our situation, $r=2$.
\end{proof}

This characterization of minimal generators agrees with what we have in Example~\ref{e.211112second}, in which the dual of the polarization of the ideal of the curve has 10 minimal generators. Note that the generators of $J^{\vee}$ are never of the form $a_ia_j$, $b_ib_j$, $c_ic_j$, or $d_id_j$. Therefore in the complementary graph $\bar{G}$, all of these edges are present. The next result gives some information about cycles in $\bar{G}$. 

\begin{lemma} \label{l.describecycle}
Let $C=(p_1,\dots,p_6)$ be a tetrahedral curve with complementary graph $\bar{G}$. Suppose there is an induced $r$-cycle in $\bar{G}$ with $r \ge 4$ (that is, there is a cycle in $\bar{G}$ of length $r \ge 4$ with no chord). Then $r=4$, and the cycle has vertices $a_i$, $b_j$, $c_l$, and $d_m$ for some $i,j,l$, and $m$.
\end{lemma}

\begin{proof}
Pick any induced $r$-cycle $C_r$ in $\bar{G}$ with $r \ge 4$. Suppose we have two vertices of one type in the cycle; that is, without loss of generality, say we have both $a_{i_1}$ and $a_{i_2}$ in the cycle with $i_1 < i_2$. Because no generator of the edge ideal of $G$ has the form $a_{i_1}a_{i_2}$, these vertices must be adjacent in $\bar{G}$, so they are adjacent in the induced cycle $C_r$. Suppose that, without loss of generality, some $b_j$, $a_{i_1}$, and $a_{i_2}$, are in a row in that order in $C_r$. Then since there is no edge between $b_j$ and $a_{i_2}$ in $\bar{G}$, there must be such an edge in $G$, meaning $a_{i_2}b_j$ is a generator of $J^{\vee}$, and $i_2+j \le p_1+1$. Moreover, since $b_j$ and $a_{i_1}$ are joined by an edge in $\bar{G}$, we conclude that $a_{i_1}b_j$ is not a generator of $J^{\vee}$. But $i_1 < i_2$, so $i_1 + j < i_2+j \le p_1+1$, which means that $a_{i_1}b_j$ is a generator of $J^{\vee}$ by the criterion of Lemma~\ref{l.dualgens}, a contradiction. Hence any induced cycle has at most one vertex of each of the types $a$, $b$, $c$, and $d$, which also implies that there are no induced $r$-cycles in $\bar{G}$ for $r \ge 5$.
\end{proof}

\begin{remark} \label{r.threevar}
Let $S=k[a,b,c]$, and let $I=(a,b)^{p_1} \cap (a,c)^{p_2} \cap (b,c)^{p_3}$ be an ideal in $S$. As in the tetrahedral curves case, we can polarize $I$ and take the Alexander dual to get a new ideal $J^{\vee}$, which has a graph $G$ and a complementary graph $\bar{G}$ associated to it. It follows from an argument similar to the one in the proof of Lemma~\ref{l.describecycle} that $\bar{G}$ has no induced $r$-cycles for $r \ge 4$ since there are only $a$, $b$, and $c$ variables in the polarization of $I$. Hence $\bar{G}$ is chordal, and $S/I$ is Cohen-Macaulay. (Of course, there are a number of other easy ways to prove this as well, including thinking of $I$ as the ideal of three general fat points in $\mathbb P^2$.)
\end{remark}

As a consequence of Lemma~\ref{l.describecycle}, to determine whether $\bar{G}$ is chordal, we need only determine whether there is an induced 4-cycle, which we know will consist of vertices $a_i$, $b_j$, $c_l$, and $d_m$ for some positive integers $i,j,l$, and $m$. 

\begin{lemma} \label{l.whichcycles}
Under the standing assumption that $p_1+p_6 \ge \max(p_2+p_5,p_3+p_4)$, $\bar{G}$ is chordal if and only if it contains no induced cycle of the form $a_ic_lb_jd_ma_i$.
\end{lemma}

\begin{proof}
One direction is trivial. For the other, suppose there is no induced cycle of the form $a_ic_lb_jd_ma_i$. We need to show that there is no other type of induced 4-cycle. The existence of a 4-cycle $a_ic_lb_jd_ma_i$ in $\bar{G}$ would correspond to having generators $a_ib_j$ and $c_ld_m$ of $J^{\vee}$ and $a_ic_l$, $a_id_m$, $b_jc_l$, and $b_jd_m \not \in J^{\vee}$. Combining these remarks with the criterion of Lemma~\ref{l.dualgens} giving the form of the generators of $J^{\vee}$, there is a $a_ic_lb_jd_ma_i$ 4-cycle that is an induced subgraph of $\bar{G}$ if and only if: \[ \begin{array}{ccc}
i+j &\le & p_1+1\\
l+m & \le & p_6+1\\
i+l & > & p_2+1\\
i+m & > & p_3+1\\
j+l & > & p_4+1\\
j+m & > & p_5+1
\end{array}
\]

The first two inequalities force $a_ib_j$ and $c_ld_m$ to be generators of $J^{\vee}$. The last four preclude $a_ic_l$, $a_id_m$, $b_jc_l$, and $b_jd_m$ from being generators of $J^{\vee}$. Therefore there is no $a_ic_lb_jd_ma_i$ 4-cycle in $\bar{G}$ if and only if there do not exist positive integers $i$, $j$, $l$, and $m$ such that $i+j = p_1+1$, $l+m = p_6+1$, $i+l \ge p_2+2$, $i+m \ge p_3+2$, $j+l \ge p_4+2$, and $j+m \ge p_5+2$. (We have put equal signs in the first two statements because if one can satisfy the last four inequalities with $i+j \le p_1+1$ and $l+m \le p_6+1$, one can do it with equality. In the four inequalities, we are using the fact that $i$, $j$, $l$, $m$, and $p_1, \dots, p_6$ are integers to convert $> p_i+ 1$ to $\ge p_i+2$.)

To prove that the lack of an induced cycle $a_ic_lb_jd_ma_i$ in $\bar{G}$ precludes the existence of any other type of induced 4-cycle, we show that there is no induced 4-cycle of type $a_ib_jc_ld_ma_i$ in $\bar{G}$, and the other case is analogous. If such an induced cycle existed, we would have $a_ic_l$ and $b_jd_m \in J^{\vee}$; moreover, $a_ib_j$, $a_id_m$, $b_jc_l$, $c_ld_m \not \in J^{\vee}$. Therefore we could find positive integers $i,j,l$, and $m$ such that: \[ \begin{array}{ccc}
i+l & = & p_2+1\\
j+m & = & p_5+1\\
i+j & \ge & p_1+2\\
i+m & \ge & p_3+2\\
j+l & \ge & p_4+2\\
l+m & \ge & p_6+2
\end{array}
\]

Hence $i+j+l+m=p_2+p_5+2$. Summing the inequalities for $p_1$ and $p_6$, we have $i+j+l+m \ge p_1+p_6+4$. Therefore $p_2+p_5+2 \ge p_1+p_6+4$, which is a contradiction since $p_2+p_5 \le p_1+p_6$.
\end{proof}

The conditions in the preceding proof for the existence of an induced cycle of the form $a_ic_lb_jd_ma_i$ lead us to the final reformulation of the question of when a tetrahedral curve is ACM.

\begin{corollary} \label{c.main}
Assuming $p_1+p_6 \ge \max(p_2+p_5,p_3+p_4)$, the tetrahedral curve $C=(p_1,\dots,p_6)$ is ACM if and only if there do not exist positive integers $i$, $j$, $l$, and $m$ such that $i+j = p_1+1$, $l+m = p_6+1$, $i+l \ge p_2+2$, $i+m \ge p_3+2$, $j+l \ge p_4+2$, and $j+m \ge p_5+2$.
\end{corollary}

\section{Some sufficient conditions for being ACM} \label{s.acm}

The problem of determining when a tetrahedral curve is ACM is now purely a numerical question about partitioning integers. In this section, we use Corollary~\ref{c.main} to identify some ACM tetrahedral curves. Throughout, we assume $p_1+p_6 \ge \max(p_2+p_5,p_3+p_4)$.

The easiest case is when either $p_1$ or $p_6$ is zero. If, for example, $p_1=0$, then in a decomposition as in Corollary~\ref{c.main}, $i+j=1$, which is impossible if both $i$ and $j$ are positive integers. Hence we can conclude:

\begin{proposition} \label{p.zero}
Let $C=(p_1,\dots,p_6)$. If $p_1$ or $p_6$ is zero, then $C$ is ACM.
\end{proposition}

The other easy case is when the difference between $p_1+p_6$ and the other two sums of that form is not large enough.

\begin{proposition} \label{p.atleasttwo}
Let $C=(p_1,\dots,p_6)$, and suppose $p_1+p_6=\epsilon + \max(p_2+p_5,p_3+p_4)$, where $\epsilon=0$ or 1. Then $C$ is ACM.
\end{proposition}

\begin{proof}
Without loss of generality, suppose $p_2+p_5 \ge p_3+p_4$ so that $p_1+p_6 = \epsilon + p_2+p_5$, where $\epsilon=0$ or 1. Suppose $C$ is not ACM, so the decomposition of Corollary~\ref{c.main} exists. Then $i+j+l+m=p_1+p_6+2$, and, summing the inequalities for $p_2$ and $p_5$, $i+j+l+m \ge p_2+p_5+4$. Hence $p_1+p_6 \ge 2+p_2+p_5$, contradicting the assumption that $p_1+p_6=\epsilon+p_2+p_5$ for $\epsilon \in \{0,1\}$. Thus $C$ is ACM.
\end{proof}

As a consequence of these results, we can classify the ACM Schwartau tetrahedral curves. Our result is stated differently from the classifications of Schwartau in \cite[Theorem 2.4]{MN} and Migliore and Nagel in \cite[Theorem 5.3]{MN} because we are assuming that $p_1+p_6$ is at least as large as $p_2+p_5$ and $p_3+p_4$. Of course, Corollary~\ref{c.schwartau} also follows from our main result, Theorem~\ref{t.main}, but it seems easier to prove it directly, and thus we do so here. The proof also serves as an introduction to the method of proof of Proposition~\ref{p.3ormore}.

\begin{corollary} \label{c.schwartau}
Let $C$ be a Schwartau tetrahedral curve $(p_1,0,p_3,p_4,0,p_6)$, and suppose without loss of generality that $p_1+p_6 \ge p_3+p_4$. Then $C$ is ACM if and only if $p_1=0$, $p_6=0$, or $p_1+p_6=\epsilon+p_3+p_4$, where $\epsilon \in \{0,1\}$.
\end{corollary}

\begin{proof}
That these cases produce an ACM tetrahedral curve is immediate from Propositions~\ref{p.zero} and ~\ref{p.atleasttwo}. We show directly that in all other cases, $C$ fails to be ACM. Assume that $p_1, p_6 > 0$ and $p_1+p_6 \ge 2+p_3+p_4$. Suppose $i$, $j$, $l$, and $m$ are positive integers such that $i+j=p_1+1$, and $l+m = p_6+1$. If $i+m \ge p_3+2$ and $j+l \ge p_4+2$, we are done, for a Corollary~\ref{c.main} decomposition exists. Otherwise, at least one of those two inequalities fails. If both fail, then $p_1+p_6+2=i+j+l+m < p_3+p_4+4$, contradicting the fact that $p_1+p_6 \ge 2+p_3+p_4$. We reach a similar contradiction if $i+m \le p_3+2$ and $j+l < p_4+2$ or $i+m < p_3+2$ and $j+l \le p_4+2$.

Therefore, without loss of generality, we may assume that $i+m \ge p_3+3$ and $j+l \le p_4+1$. As long as $i$ and $m$ are not both 1, we may either decrease $i$ by one and increase $j$ by one or decrease $m$ by one and increase $l$ by one. This maintains the partitioning of $p_1+1$ and $p_6+1$, and, by induction, we may repeat this process as long as necessary to ensure that $i+m \ge p_3+2$ and $j+l \ge p_4+2$. If $i=1=m$, then we are stuck because we cannot shift from $i$ or $m$ to $j$ or $l$. But in that case, $2=i+m \ge p_3+3$, so $-1 \ge p_3$, a contradiction since $p_3$ is nonnegative. Hence a Corollary~\ref{c.main} decomposition exists, and $C$ fails to be ACM.
\end{proof}

There is one more situation in which we can easily conclude that $C$ is ACM. The conditions in the next proposition are not particularly tidy, but the proof is easy, and isolating these inequalities is vital in classifying the final sufficient condition for a tetrahedral curve to be ACM, which we prove in the next section.

\begin{proposition} \label{p.eachletter}
Let $C=(p_1,\dots,p_6)$ be a tetrahedral curve with $p_1+p_6$ maximal among $p_1+p_6$, $p_2+p_5$, and $p_3+p_4$. If the $p_i$ fail to satisfy any one of the inequalities below, then $C$ is ACM. \[ \begin{array}{ccc}
2p_1 & \ge & p_2+p_3+3-p_6\\
2p_1 & \ge & p_4+p_5+3-p_6\\
2p_6 & \ge & p_2+p_4+3-p_1\\
2p_6 & \ge & p_3+p_5+3-p_1
\end{array}\]
\end{proposition}

\begin{proof}
We prove that not satisfying the first inequality implies that $C$ is ACM, and the others are similar. The conditions on $i$, $j$, $l$, and $m$ from Corollary~\ref{c.main} give that $i+j=p_1+1$, $i+l \ge p_2+2$, and $i+m \ge p_3+2$. Therefore $3i+j+l+m \ge p_1+p_2+p_3+5$. Using the fact that $i+j+l+m=p_1+p_6+2$, we have $2i+p_1+p_6+2 \ge p_1+p_2+p_3+5$, and therefore $2i \ge p_2+p_3+3-p_6.$ If a Corollary~\ref{c.main} decomposition exists, then $p_1 \ge i$, so we conclude that if $C$ is not ACM, then $2p_1 \ge p_2+p_3+3-p_6$. Hence if that inequality fails, $C$ is ACM.
\end{proof}

\begin{example} \label{e.eachletter}
The tetrahedral curve $(1,1,1,3,2,5)$ is an ACM curve. Note that \[ 2p_1 = 2 < 3=3+2+3-5 = p_4+p_5+3-p_6,\] so the second inequality of Proposition~\ref{p.eachletter} fails. (In fact, $(1,0,0,3,2,5)$ is still ACM since lowering $p_2$ and $p_3$ makes no difference.)
\end{example}

\section{Final cases} \label{s.final}

In our final section, we find one more condition under which a tetrahedral curve is ACM and give the complete list of necessary and sufficient conditions. We begin with a proposition that significantly reduces the number of remaining curves to consider, identifying a collection of tetrahedral curves that are not ACM. 

\begin{proposition} \label{p.3ormore}
Let $C=(p_1,\dots,p_6)$ be a tetrahedral curve such that $p_1+p_6 \ge 2+\max(p_2+p_5,p_3+p_4)$ and $p_1, p_6 > 0$. Suppose further that $p_1+p_6 \ge 3+ \min(p_2+p_5,p_3+p_4)$, and $C$ satisfies all the inequalities of Proposition~\ref{p.eachletter}. Then $C$ is not ACM.
\end{proposition}

\begin{proof}
Let $i,j,l,$ and $m$ be positive integers such that $i+j=p_1+1$ and $l+m=p_6+1$. We may assume that $p_2+p_5 \le p_3+p_4$, and thus $p_1+p_6 \ge 3+p_2+p_5$. We claim first that by an argument identical to that used in Corollary~\ref{c.schwartau}, we may choose $i,j,l$, and $m$ such that $i+l \ge p_2+2$ and $j+m \ge p_5+2$. To see this, suppose instead that $i+l \le p_2+2$ and $j+m < p_5+2$. Then $p_1+p_6+2 = i+j+l+m < p_2+p_5+4$, so $p_1+p_6 < p_2+p_5+2$, a contradiction; similarly, $i+l < p_2+2$ and $j+m \le p_5+2$ is impossible. Thus, if we do not have both $i+l \ge p_2+2$ and $j+m \ge p_5+2$, we must have $i+l \ge p_2+3$ and $j+m \le p_5+1$ or $i+l \le p_2+1$ and $j+m \ge p_5+3$. Without loss of generality, assume the former. In that case, unless $i=l=1$, we can lower $i$ and raise $j$ or lower $l$ and raise $m$ until we get $i+l \ge p_2+2$ and $j+m \ge p_5+2$. This does not affect the equalities $i+j=p_1+1$ and $l+m=p_6+1$. If $i=l=1$, then $2=i+l \ge p_2+3$, which would say that $-1 \ge p_2$, a contradiction.

Therefore we may choose $i$, $j$, $l$, and $m$ such that $i+j=p_1+1$, $l+m=p_6+1$, $i+l \ge p_2+2$, and $j+m \ge p_5+2$. Suppose the two inequalities were both equalities. Then $p_1+p_6+2=i+j+l+m = p_2+p_5+4$, contradicting the assumption that $p_1+p_6 \ge 3+p_2+p_5$. Hence, without loss of generality, we may choose $i,j,l$, and $m$ such that $i+j=p_1+1$, $l+m=p_6+1$, $i+l \ge p_2+3$, and $j+m \ge p_5+2$. Now we need only ensure that the inequalities $i+m \ge p_3+2$ and $j+l \ge p_4+2$ are satisfied.

If those inequalities already hold under our choices of $i$, $j$, $l$, and $m$, we are done. Suppose not. One possibility is that $i+m \le p_3+2$ and $j+l < p_4+2$. In this case, $p_1+p_6+2 = i+j+l+m < p_3+p_4+4$, so $p_1+p_6 < p_3+p_4+2$, a contradiction. The analogous argument rules out $i+m < p_3+2$ and $j+l \le p_4+2$.

The remaining possibilities are $i+m \ge p_3+3$ and $j+l \le p_4+1$ or $i+m \le p_3+1$ and $j+l \ge p_4+3$; we assume the former. Our goal is to lower $i$ and raise $j$ and/or lower $m$ and raise $l$ until we have $i+m \ge p_3+2$ and $j+l \ge p_4+2$. None of these actions affect the equations $i+j=p_1+1$ and $l+m=p_6+1$. The only impediments to these moves are if $i=1=m$, or if decreasing $i$ or $m$ would cause one of the inequalities $i+l \ge p_2+2$ or $j+m \ge p_5+2$ to fail. As long as $i >1$, we can decrease $i$ by one and increase $j$ by one since $i+m \ge p_3+3$ and $i+l \ge p_2+3$, meaning each of those inequalities has some leeway. 

Suppose now instead that we reach the point at which $i=1$. If also $m=1$, then $2=i+m \ge p_3+3$, so $-1 \ge p_3$, a contradiction. Assume instead that $i=1$ but $m \not = 1$. Because $i+j=p_1+1$ and $i=1$, we know $j=p_1$. Since $m \not = 1$, we can lower $m$ and increase $l$ unless $j+m=p_5+2$. If $j+m=p_5+2$, we have $m=p_5+2-p_1$, and $l=p_6+1-m=p_6+1-(p_5+2-p_1)=p_6-p_5+p_1-1$. The inequality $j+l \le p_4+1$ tells us that \[ j+l = p_1 + (p_6-p_5+p_1-1) = 2p_1 + p_6-p_5-1 \le p_4+1,\] meaning that \[ 2p_1 \le p_4+p_5+2-p_6.\] But we are assuming that the inequalities of Proposition~\ref{p.eachletter} all hold, so in particular, $2p_1 \ge p_4+p_5+3-p_6$, a contradiction.

If instead $i+m \le p_3+1$ and $j+l \ge p_4+3$, the argument is virtually the same, this time using the inequality $2p_6 \ge p_3+p_5+3-p_1$ from Proposition~\ref{p.eachletter} at the end.
\end{proof}

In view of Propositions~\ref{p.atleasttwo} and ~\ref{p.3ormore}, we now need only investigate what happens when $p_1+p_6 = 2+p_2+p_5 = 2+p_3+p_4$. The next proposition illustrates the remaining case in which $C$ is ACM. We may assume that $p_1, p_6 > 0$ and that the inequalities of Proposition~\ref{p.eachletter} all hold because otherwise, we know $C$ is ACM.

\begin{proposition} \label{p.lastcase}
Let $C=(p_1,\dots,p_6)$ be a tetrahedral curve with $p_1,p_6 > 0$ and $p_1+p_6 = 2+p_2+p_5=2+p_3+p_4$. Suppose all the inequalities of Proposition~\ref{p.eachletter} hold. Then $C$ is ACM if and only if $p_1+p_3+p_5$ is even.
\end{proposition}

\begin{proof}
As before, we may choose positive integers $i,j,l$, and $m$ such that $i+j=p_1+1$, $l+m=p_6+1$, $i+l \ge p_2+2$ and $j+m \ge p_5+2$. Therefore \[ p_1+p_6+2 = i+j+l+m \ge p_2+p_5+4,\] meaning that $p_1+p_6 \ge 2+p_2+p_5$. But we know equality holds, and therefore $i+l = p_2+2$ and $j+m = p_5+2$. If a decomposition as in Corollary~\ref{c.main} exists, then by a similar argument, we must have $i+m=p_3+2$ and $j+l=p_4+2$. Therefore $C$ has a Corollary~\ref{c.main} decomposition if and only if there exist positive integers $i,j,l$, and $m$ such that: \[
\begin{array}{ccc}
i+j & = & p_1+1\\
i+l & = & p_2+2\\
i+m & = & p_3+2\\
j+l & = & p_4+2\\
j+m & = & p_5+2\\
l+m & = & p_6+1
\end{array}
\]
Doing Gaussian elimination on these six equations in four variables yields the following matrix in reduced row-echelon form: \[\begin{pmatrix}1 & 0 & 0 & 0 & \vline & \frac{1}{2}(p_1+p_3-p_5+1)\\ 0 & 1 & 0 & 0 & \vline & \frac{1}{2}(p_1-p_3+p_5+1) \\ 0 & 0 & 1 & 0 & \vline & \frac{1}{2}(-p_1+2p_2-p_3+p_5+3)\\ 0 & 0 & 0 & 1 & \vline & \frac{1}{2}(-p_1+p_3+p_5+3) \\ 0 & 0 & 0 & 0 & \vline & p_1-p_3-p_4+p_6-2 \\ 0 & 0 & 0 & 0 & \vline & p_1-p_2-p_5+p_6-2 \end{pmatrix}
\]

Of course, there are many equivalent formulations of the expressions on the right-hand side due to the relations $p_1+p_6=2+p_2+p_5=2+p_3+p_4$, and the last two rows just restate those relations. The other four rows give us the unique solution for $i,j,l$, and $m$: \[\begin{array}{ccc}
i & = & \frac{1}{2}(p_1+p_3-p_5+1)\\
j & = & \frac{1}{2}(p_1-p_3+p_5+1)\\
l & = & \frac{1}{2}(-p_1+2p_2-p_3+p_5+3)\\
m & = & \frac{1}{2}(-p_1+p_3+p_5+3)
\end{array}
\]

This is a decomposition in the sense of Corollary~\ref{c.main}, and hence $C$ is not ACM, if and only if each formula for $i,j,l$, and $m$ yields a positive integer. The expressions are clearly all integers if and only if $p_1+p_3+p_5$ is odd. We show here that the formulas for $j$ and $l$ are positive, and the other two cases are analogous.

To show that the formula for $j$ is positive, we prove that $p_1-p_3+p_5+1 \ge 2$; this is stronger than what we need, but if a Corollary~\ref{c.main} decomposition exists, $j$ would be an integer, and thus this inequality would hold. We are assuming that the inequalities of Proposition~\ref{p.eachletter} are satisfied, so $2p_1+p_6 \ge p_2+p_3+3$. Using the fact that $p_1+p_6=2+p_2+p_5$, we have \[ p_1 + 2+ p_2+p_5 \ge p_2+p_3+3,\] and hence \[ p_1+p_5 \ge p_3+1,\] which is equivalent to the statement that $p_1-p_3+p_5+1 \ge 2$.

Finally, we prove that the formula for $l$ is always positive, showing that $-p_1+2p_2-p_3+p_5+3 \ge 2$, or $p_1+p_3 \le 2p_2+p_5 +1$. Because the inequalities of Proposition~\ref{p.eachletter} hold, $2p_6 + p_1 \ge p_3+p_5+3$. Therefore, using $p_1+p_6=2+p_2+p_5$, \[ p_6 +(p_1+p_6-2) \ge p_3+p_5+3-2 \iff p_6+p_2+p_5 \ge p_3+p_5+1,\] so $p_6+p_2 \ge p_3+1$. Adding $p_1$ to both sides, we get \[ p_1+p_6+p_2 \ge p_1+p_3+1 \iff p_2+p_5+2+p_2 \ge p_1+p_3+1,\] and hence \[p_1+p_3 \le 2p_2+p_5+1.\]
\end{proof}

We are now ready to prove our main theorem.

\begin{theorem} \label{t.main}
Let $C$ be a tetrahedral curve with exponent vector $(p_1,\dots,p_6)$. Suppose without loss of generality that $p_1+p_6=\max(p_1+p_6,p_2+p_5,p_3+p_4)$. Then $C$ is ACM if and only if at least one of the following conditions holds:
\begin{enumerate}[(i)]
\item $p_1=0$ or $p_6=0$ 
\item $p_1+p_6=\epsilon+\max(p_2+p_5,p_3+p_4)$, where $\epsilon \in \{0,1\}$
\item $2p_1 < p_2+p_3+3-p_6$ or $2p_1 < p_4+p_5+3-p_6$ or $2p_6 < p_2+p_4+3-p_1$ or $2p_6 < p_3+p_5+3-p_1$
\item All inequalities of (iii) fail, $p_1+p_6=2+p_2+p_5=2+p_3+p_4$, and $p_1+p_3+p_5$ is even.
\end{enumerate}
\end{theorem}

\begin{proof}
The cases listed give ACM curves by Propositions~\ref{p.zero}, ~\ref{p.atleasttwo}, ~\ref{p.eachletter}, and ~\ref{p.lastcase}. If we are not in any of those cases, then either $p_1+p_6 \ge 3+\max(p_2+p_5,p_3+p_4)$, and it follows from Proposition~\ref{p.3ormore} that $C$ is not ACM, or $p_1+p_6=2+\max(p_2+p_5,p_3+p_5)$. In that case, if $p_1+p_6 \ge 3+\min(p_2+p_5,p_3+p_4)$, $C$ fails to be ACM by Proposition~\ref{p.3ormore}; otherwise, $p_1+p_6=2+p_2+p_5=2+p_3+p_4$. Then Proposition~\ref{p.lastcase} proves that $C$ is ACM if and only if $p_1+p_3+p_5$ is even.
\end{proof}

We hope that the illustration of this combination of combinatorial methods in this paper will help stimulate ongoing research on resolutions of facet ideals and attempts to generalize Fr\"oberg's Theorem~\ref{t.froberg}. Successes in these areas would help yield future results in the spirit of this work for more general monomial ideals. 



\begin{thebibliography}{99}

\bibitem{BH} W.~Bruns and J.~Herzog, \emph{Cohen-{M}acaulay rings}, in \emph{Cambridge Studies in Advanced Mathematics} \textbf{39} (Cambridge University Press, Cambridge, 1993).

\bibitem{ER} J. Eagon and V. Reiner, Resolutions of Stanley-Reisner rings and Alexander duality. \emph{J. Pure Appl. Algebra} {\bf 130} (1998), no. 3, 265--275.
\bibitem{Faridi-facet} S. Faridi, The facet ideal of a simplicial complex. \emph{Manuscripta Math.} {\bf 109} (2002), 159--174.

\bibitem{Faridi} S. Faridi, Monomial ideals via square-free monomial ideals. \emph{Commutative algebra}, Lect. Notes Pure Appl. Math. {\bf 244}, Chapman \& Hall/CRC, Boca Raton, FL, 2006, 85--114.

\bibitem{FMN} C. Francisco, J. Migliore, and U. Nagel, On the componentwise linearity and the minimal free resolution of a tetrahedral curve. \emph{J. Algebra} {\bf 299} (2006), no. 2, 535--569. 

\bibitem{Froberg} R. Fr\"oberg, On Stanley-Reisner rings. In Topics in algebra, Banach Center Publications, {\bf 26} (2) (1990), 57--70.

\bibitem{M2}
D.~R. Grayson and M.~E. Stillman, \emph{Macaulay 2, a software system for
  research in algebraic geometry}.
\newblock \verb|http://www.math.uiuc.edu/Macaulay2/|.

\bibitem{HVT} H. T\`ai H\`a, A. Van Tuyl, Splittable ideals and the resolutions
of monomial ideals. \emph{J. Algebra} {\bf 309} (2007), 405--425.

\bibitem{HVTsurvey} H. T\`ai H\`a, A. Van Tuyl, Resolutions of square-free monomial ideals via facet ideals: a survey. Preprint, 2006. {\tt math.AC/0604301}

\bibitem{HHZ} J. Herzog, T. Hibi, and X. Zheng, Dirac's theorem on chordal graphs and Alexander duality. \emph{European J. Combin.} {\bf 25} (2004), no. 7, 949--960.

\bibitem{Katzman} M. Katzman, Characteristic-independence of Betti numbers of graph ideals. \emph{J. Combin. Theory Ser. A} {\bf 113} (2006), no. 3, 434--454.

\bibitem{MN} J. Migliore and U. Nagel, Tetrahedral curves. 
\emph{Int. Math. Res. Notices} {\bf 15} (2005), 899--939.

\bibitem{S} P.\ Schwartau, {\em Liaison Addition and Monomial
Ideals},  Ph.D.\ thesis, Brandeis University (1982).

\bibitem{V} R. Villarreal, \emph{Monomial algebras}, in \emph{Monographs and Textbooks in Pure and Applied Mathematics} {\bf 238} (Marcel Dekker, Inc., New York, 2001)

\bibitem{Zheng} X. Zheng, Resolutions of facet ideals. \emph{Comm. Algebra} {\bf 32} (2004), no. 6, 2301--2324.
\end{thebibliography}
\end{document}